\documentclass{article}
\usepackage[utf8]{inputenc}
\usepackage[english]{babel}
\usepackage{amsmath,amsfonts} 
\usepackage{amsthm} 
\usepackage{amssymb}
\usepackage{graphicx}
\usepackage{paralist}
\usepackage{mathtools}
\usepackage[bb=boondox]{mathalfa}

\usepackage{xcolor}
\usepackage{subfig}
\definecolor{ForestGreen}{RGB}{34,139,34}

\theoremstyle{definition}
\newtheorem{definition}{Definition}		
\newtheorem{example}[definition]{Example}

\theoremstyle{plain}
\newtheorem{lemma}[definition]{Lemma}				
\newtheorem{corollary}[definition]{Corollary}		
\newtheorem{theorem}[definition]{Theorem}

\theoremstyle{remark}

\numberwithin{equation}{section}

\newcommand{\R}{\mathbb{R}} 
\newcommand{\N}{\mathbb{N}} 

\DeclareMathOperator{\id}{Id}

\usepackage[affil-it]{authblk}

\title{A measure--theoretic representation of graphs}

\author[1,2]{Raffaella Mulas}
\author[1]{Giulio Zucal\thanks{zucal@mis.mpg.de}}
\affil[1]{Max Planck Institute for Mathematics in the Sciences, Leipzig, Germany}
\affil[2]{Vrije Universiteit Amsterdam, Amsterdam, The Netherlands}

\date{}

\begin{document}

\maketitle
\begin{abstract}
Inspired by the notion of action convergence in graph limit theory, we introduce a measure--theoretic
representation of matrices, and we use it to define a new notion of pseudo-metric on the space of matrices. Moreover, we show that such pseudo-metric is a metric on the subspace of adjacency or Laplacian 
matrices for graphs. Hence, in particular, we obtain a metric for isomorphism classes of graphs. Additionally, we study how some properties of graphs translate in this measure representation, and we show how our analysis contributes to a simpler understanding of action convergence of graphops.

\vspace{0.2cm}
\noindent {\bf Keywords:} 
Network distances, Graph limits, Action convergence, Graphops
\end{abstract}

\section{Introduction}
 In the last 20 years, the study of complex networks has permeated many areas of social and natural sciences. Important examples are computer, telecommunication, biological, cognitive, semantic and social networks. 
Network structures are usually modeled using graph theory to represent
pairwise interactions between elements of a network. For this reason, it is particularly important to find advantageous ways to represent and compare graphs. Many graph distances have been proposed from both an applied and a theoretical perspective.  In applications, the most common pseudo-distances are inspired by local comparison (e.g.\
Hamming distance, Jaccard distance) and/or global spectral methods. For an overview of the most commonly-used graph pseudo-distances to compare empirical networks in practice, see  \cite{NetworkDistancesApplied}. From a theoretical point of view, the selection of a metric on the space of graphs is related to graph limit theory. This is a very active field of mathematics that connects graph theory with many other mathematical areas such as stochastic processes, ergodic theory, spectral theory and several branches of analysis and topology. In fact, in mathematical terms, one is interested in finding a metric/topology on the space of graphs, and a completion of this space with respect to that topology. Traditionally, this field grew in two distinct directions: limits of graph sequences of bounded degree on the one hand (Benjamini-Schramm convergence \cite{BenjaminiLimit} and the stronger notion of local-global convergence \cite{local-global1,Hatami2014LimitsOL}), and limits of dense graph sequences on the other hand (graphons and the related notion of cut-metric \cite{BORGS20081801,Lovsz2007SzemerdisLF, LOVASZ2006933}). For a complete treatment of these topics, see the monograph by L.\ Lovász \cite{LovaszGraphLimits}. More recently, the challenging intermediate case of sparse graph sequences with unbounded degree attracted a lot of interest, as this case covers the vast majority of networks in applications. In fact, real networks are usually sparse, but not very sparse and heterogeneous. For instance, there is the relaxation of graphons to less dense graph sequences, namely $L^p$--graphons \cite{LpGraphons1,Lpgraphon2}. In a recent paper, A.\ Backhausz and B.\ Szegedy introduced a new functional analytic/measure--theoretic notion of convergence \cite{backhausz2018action}, that not only covers the intermediate degree case, but also unifies the graph limit theories that we mentioned previously. Other works in this direction are \cite{10.5555/3122009.3242067,Borgs2019,KUNSZENTIKOVACS2022109284,MarkovSpaces,VeitchVictor}.\newline

In this paper, we contribute to the study of representation and comparison of graphs, by introducing and investigating the following measure--theoretic representation of matrices. \begin{definition}\label{defMeasGenerated}
Let $A$ be an $n\times n$ matrix and, given $\mathbf{x}=(x_i)_i\in \mathbb{R}^n$, let
\begin{equation*}
    \mu^{\mathbf{p}}_{(A,\,\mathbf{x})}\coloneqq \sum_{i=1}^n p_i\cdot \delta_{(x_i,(A\mathbf{x})_i)},
\end{equation*}where $\delta$ denotes the Dirac measure, and $p_1,\ldots,p_n\in\mathbb{R}_{>0}$ such that $\sum_{i=1}^n p_i=1$ are fixed. Let also $\mathbf{p}\coloneqq (p_i)_i\in\mathbb{R}^n$. We will simply denote by $\mu_{(A,\,\mathbf{x})}$ the measure $\mu^{\mathbf{p}}_{(A,\,\mathbf{x})}$, where there is no risk of confusion about the given probability vector $\mathbf{p}$.\newline
The \emph{family of measures generated by $A$} is $$(\mu_{(A,\,\mathbf{x})})_{\mathbf{x}\in \mathbb{R}^n }.$$
The \emph{set of measures generated by $A$} is
$$\mathcal{Z}_A\coloneqq \{\mu_{(A,\,\mathbf{x})}:\mathbf{x}\in \mathbb{R}^n\}.$$We will drop the $A$ in the subscript of $\mu_{(A,\,\mathbf{x})}$ and $\mathcal{Z}_A$, writing  $\mu_\mathbf{x}$ and $\mathcal{Z}$, respectively, whenever the dependence on the matrix $A$ is obvious.
\end{definition}

We use this representation to define a new notion of pseudo-metric on the space of matrices. Moreover, we show that, for this pseudo-metric, the distance  between two matrices $A$ and $B$ is zero if and only if $A$ and $B$ are \emph{switching equivalent}, i.e., there exists a permutation matrix $P$ such that $A=PBP^\top$. Formally,

\begin{theorem}\label{thm:main}
Let $K\in \mathbb{R}_{\geq 1}$ and let $\mathcal{S}$ be the set of $n\times n$ matrices $A$ such that $||A||_{\infty \rightarrow 1}\leq K$. Then, a matrix $A\in \mathcal{S}$ is determined, up to switching equivalence, by the set $\mathcal{Z}_A$ of measures generated by $A$, where we consider the measures relative to the uniform probability measure.
\end{theorem}

For example, if two matrices $A$ and $B$ are the adjacency (or Laplacian) matrices of two graphs $G_1$ and $G_2$, respectively, then they are switching equivalent if and only if $G_1$ and $G_2$ are isomorphic. As a consequence, the above framework will allow us to define a metric on the class of isomorphic graphs.\newline
Additionally, we study how some properties of graphs, as the spectrum, the vertex degrees, and some homomorphism numbers, translate in this measure representation.\newline 

Such representation and metric are inspired by the notion of action convergence in graph limits theory \cite{backhausz2018action}, but it is simpler. For this reason, our analysis also contributes to a simpler understanding of the limit notion for graphs. Our results show that for discrete probability measures with finite support, such distance has the same expressive power of the more complex action convergence metric \cite[Definition 2.6]{backhausz2018action}. However, for general measures it stays open whether this metric has strictly less expressive power or not. 

\subsection*{Structure of the paper}
This work is structured as follows. In Section \ref{Definitions} we introduce some relevant definitions and notations, in Section \ref{PrelResults} we give some preliminary results, and in Section \ref{MaResults} we prove Theorem \ref{thm:main}. Moreover, in Section \ref{PropertiesGraphs} we relate properties of matrices and graphs with the measure representation, and in Section \ref{RelatActionConv} we underline the relationship with action convergence.  Finally, in Section \ref{Concl}, we present some open questions as well as future directions.

\section{Basic definitions and notations}

\label{Definitions}Throughout the paper we fix $n\in \mathbb{N}_{\geq 2}$, we let $\mathcal{P}$ denote the set of permutation matrices of order $n$, and we let $\mathbf{e_i}$ denote the $i$-th vector of the canonical basis of $\mathbb{R}^n$, for $i=1,\ldots,n$. We use $\delta$ to denote a Dirac measure and $\mathcal{P}\left(\mathbb{R}^{p}\right)$ to represent the space of probability measures in $\mathbb{R}^p$.
\begin{definition}
Let $\rho$ be a probability measure on $\mathbb{R}^2$ whose support is given by exactly $n$ points. A pair of vectors $(\mathbf{x},\mathbf{y})$ with $\mathbf{x}=(x_i)_i, \mathbf{y}=(y_i)_i\in\mathbb{R}^n$ is an \emph{ordered support} of $\rho$ if
\begin{equation*}
    x_1\leq \ldots\leq x_n
\end{equation*}and
\begin{equation*}
    \rho=\sum_{i=1}^n p_i\cdot \delta_{(x_i,y_i)},
\end{equation*}for some $p_1,\ldots,p_n\in\mathbb{R}_{>0}$ such that $\sum_{i=1}^n p_i=1$.
\end{definition}

Observe that, in the above definition, if the entries of $\mathbf{x}$ are all pairwise distinct, then there is a unique ordered support of $\rho$.\newline

From here on, we also fix $p_1,\ldots,p_n\in\mathbb{R}_{>0}$ such that $\sum_{i=1}^n p_i=1$, and we let $\mathbf{p}\coloneqq (p_i)_i\in\mathbb{R}^n$.\newline

In the following, we will mainly focus on the case where $\mathbf{p}$ is the uniform distribution on $[n]$. However, we point out that, in some cases, other probability measures can be more appropriate. For example, if we consider the family of measures generated by the transition probability matrix of a reversible Markov chain, an appropriate choice for $\mathbf{p}$ would be the stationary distribution of the Markov chain.\newline

Given an $n\times n$ matrix $A$ and a vector $\mathbf{x}=(x_i)_i\in \mathbb{R}^n$, we
define the \emph{marginal with respect to the first variable} as the discrete measure $\mu^{\mathbf{p},M}_\mathbf{x}$ such that \begin{equation*}
\begin{aligned}
     \mu^{\mathbf{p},M}_\mathbf{x}(\{q_1\})\coloneqq \sum_{q_2\in \R} \left(\sum_{i=1}^n p_i\cdot \delta_{(x_i,(A\mathbf{x})_i)}(\{(q_1,q_2)\})\right)=\\
\sum_{q_2\in \{(A\mathbf{x})_j: \ j\in [n]\}}\left(\sum_{i=1}^n p_i\cdot \delta_{(x_i,(A\mathbf{x})_i)}(\{(q_1,q_2)\})\right)
\end{aligned}\end{equation*}for every $q_1\in \R$. Being the measure discrete, it is completely characterized by the previous formula. Also in this case, we will denote by $\mu^{M}_\mathbf{x}$ the measure $\mu^{\mathbf{p},M}_\mathbf{x}$ when there is no risk for confusion. We also let 
\begin{align*}
    \mathcal{Z}_{(A,\,\mathbf{x})}\coloneqq \mathcal{Z}_{\mathbf{x}}\coloneqq& \left\{\mu_{\mathbf{y}}\in \mathcal{Z}: \mu^M_{\mathbf{y}}=\sum_{i=1}^n p_i\cdot \delta_{x_i}\right\}\\ =&\bigl\{\mu_{P\mathbf{x}}\in\mathcal{Z}:\, P\in\mathcal{P} \text{ and } P\mathbf{p}=\mathbf{p}\bigr\}.
\end{align*}

We observe that $\mathcal{Z}_{\mathbf{x}}$ has finitely many elements, since $\mathcal{P}$ is a finite set, and that $\mathcal{Z}_{\mathbf{x}}=\mathcal{Z}_{P\mathbf{x}}$ for all $P\in \mathcal{P}$.

\begin{example}
As an easy example, assume that $n=2$ and $p_1=p_2=\frac{1}{2}$. Then, $\mathcal{P}=\{\id,P\}$, where
\begin{equation*}
    \id\coloneqq \begin{pmatrix}1 & 0\\ 0 & 1
    \end{pmatrix} \quad \text{and}\quad P\coloneqq \begin{pmatrix}0 & 1\\ 1 & 0
    \end{pmatrix}.
\end{equation*}Let also
\begin{equation*}
    A\coloneqq \begin{pmatrix}a_{11} & a_{12}\\ a_{21} & a_{22}
    \end{pmatrix}.
\end{equation*}Then,
\begin{equation*}
    \mu_{\mathbf{e_1}}=\mu_{\id\mathbf{e_1}}=\frac{1}{2}\cdot\bigl(\delta_{(1,a_{11})}+\delta_{(0,a_{21})}\bigr)
\end{equation*}and
\begin{equation*}
   \mu_{\mathbf{e_2}}= \mu_{P\mathbf{e_1}}=\frac{1}{2}\cdot\bigl(\delta_{(0,a_{12})}+\delta_{(1,a_{22})}\bigr).
\end{equation*}In particular,
\begin{equation*}
    \mathcal{Z}_{\mathbf{e_1}}=\mathcal{Z}_{\mathbf{e_2}}=\{\mu_{\mathbf{e_1}},\mu_{\mathbf{e_2}}\}.
\end{equation*}
\end{example}

The following vectors will play an important role in the proof of our main result. Therefore, we give them a name. 

\begin{definition}Let $A$ be an $n\times n$ matrix and let  $(\mu_\mathbf{x})_{\mathbf{x}\in \mathbb{R}^n }$ be the family of measures generated by $A$. A vector $\mathbf{v}\in \R^n$ is $\mathbf{p}-$\emph{irreducible for $A$} if the following condition holds. For every $P_1\in \mathcal{P}$ and for every vector $\mathbf{y}\in \R^n$,
$$
\mu^{\mathbf{p}}_{P_1\mathbf{v}}=\mu^{\mathbf{p}}_{\mathbf{y}}
$$
if and only if there exists $P_2\in \mathcal{P}$ such that $\mathbf{y}=P_2P_1\mathbf{v}$, $P_2A=AP_2$ and $P_2\mathbf{p}=\mathbf{p}$. We call it \emph{irreducible for $A$} if it is  $\mathbf{u}-$\emph{irreducible for $A$}, where $\mathbf{u}$ is the uniform n-dimensional probability vector.
\end{definition}

Notice that, if a vector is irreducible for $A$, then it is $\mathbf{p}-$irreducible for every probability vector $\mathbf{p}$. For this reason, we will just consider irreducible vectors for $A$ in our following arguments.\newline

We consider the following
\begin{example}
For the matrix 
\begin{equation*}
A=\begin{bmatrix}
0 & 2 \\
3 & 1
\end{bmatrix}, 
\end{equation*}
the vector

\begin{equation*}
\mathbf{v}=\begin{bmatrix}
1\\
0
\end{bmatrix}
\end{equation*}
is irreducible, while the vector 
\begin{equation*}
\mathbf{x}=\begin{bmatrix}
-1\\
1
\end{bmatrix}
\end{equation*}is not an irreducible vector.
\end{example}
We now want to be able to compare measures. For this reason, we recall the following well-known metric:

\begin{definition}[Lévy-Prokhorov metric]
 The \emph{Lévy-Prokhorov Metric} $d_{\mathrm{LP}}$ on the space of probability measures $\mathcal{P}\left(\mathbb{R}^{p}\right)$ is
$$\begin{aligned}
d_{\mathrm{LP}}\left(\eta_{1}, \eta_{2}\right)=&\inf \left\{\varepsilon>0: \eta_{1}(U) \leq \eta_{2}\left(U^{\varepsilon}\right)+\varepsilon \text{ and } \right.\\
&\left.\eta_{2}(U) \leq \eta_{1}\left(U^{\varepsilon}\right)+\varepsilon  \text{ for all } U \in \mathcal{B}_{p}\right\},
\end{aligned}$$

 where $\mathcal{B}_{p}$ is the Borel $\sigma$-algebra on $\mathbb{R}^{p}$ and $U^{\varepsilon}$ is the set of points that have distance smaller than $\varepsilon$ from $U$.
\end{definition}

The above metric metrizes the weak/narrow convergence for measures. \newline
Now, we want to be able to compare sets of measures. We therefore introduce the following

\begin{definition}[Hausdorff metric]
 Given $X, Y\subset \mathcal{P}\left(\mathbb{R}^{p}\right)$, their \emph{Hausdorff distance} 
$$
d_{H}(X, Y):=\max \left\{\sup _{x \in X} \inf _{y \in Y} d_{\mathrm{LP}}(x, y), \sup _{y \in Y} \inf _{x \in X} d_{\mathrm{LP}}(x, y)\right\}
$$
Note that $d_{H}(X, Y)=0$ if and only if $\operatorname{cl}(X)=\operatorname{cl}(Y)$, where $\operatorname{cl}$ is the closure in $d_{\mathrm{LP}} .$ It follows that $d_{H}$ is a pseudometric for all subsets in $\mathcal{P}\left(\mathbb{R}^{k}\right)$, and it is a metric for closed sets.
\end{definition}
By definition, the Lévy-Prokhorov distance between probability measures is upper-bounded by $1$ and, therefore, the Hausdorff metric for sets of measures is upper-bounded by one.

Now, we define the $L^\infty$ and the $L^1$ norm of a vector $\mathbf{v}\in \R^n$ as 
$$\left\|\mathbf{v}\right\|_{\infty}:=\max_{i\in [n]}|v_i|$$
and 
$$\left\|\mathbf{v}\right\|_{1}:=\sum_{i\in [n]}p_i|v_i|,$$
respectively.\newline

Additionally, we define the $(\infty \rightarrow 1)-$\emph{operator norm} of a $n\times n$ matrix $A$ as 
$$
\|A\|_{\infty \rightarrow 1}:=\sup _{ \mathrm{v}\in \R^n, \substack{ \mathrm{v}\neq 0}}\frac{\|A\mathrm{v}\|_{1} }{\|\mathrm{v}\|_{\infty}}.
$$
\section{Preliminary results }
\label{PrelResults}
In this section we prove some preliminary results that will be needed for the proof of the main one.

\begin{lemma}\label{lemma:P}
Let $A$ be an $n\times n$ matrix and let  $(\mu_\mathbf{x})_{\mathbf{x}\in \mathbb{R}^n }$ be the family of measures generated by $A$. Then, $\mathbf{x}, \mathbf{y}\in \mathbb{R}^n$ are such that
$$\mu_\mathbf{x}=\mu_\mathbf{y}$$ if and only if there exists $P\in \mathcal{P}$ such that
\begin{equation*}
    \mathbf{y}=P\mathbf{x}, \quad P\mathbf{p}=\mathbf{p} \quad \text{and }\quad \mathbf{x}\in \ker(PA-AP).
\end{equation*}
\end{lemma}
\proof
By definition, $\mu_\mathbf{x}=\mu_\mathbf{y}$ if and only if
$$\sum_{i=1}^n p_i\cdot \delta_{(x_i,(A\mathbf{x})_i)}=\sum_{i=1}^n p_i\cdot \delta_{(y_i,(A\mathbf{y})_i)},$$
hence, if and only if there exists $P\in \mathcal{P}$ such that $(P\mathbf{x},PA\mathbf{x})=(\mathbf{y},A\mathbf{y})$ and $P\mathbf{p}=\mathbf{p}$.  This happens if and only if $P\mathbf{p}=\mathbf{p}$, $\mathbf{y}=P\mathbf{x}$ and $PA\mathbf{x}=AP\mathbf{x}$, hence, if and only if $P\mathbf{p}=\mathbf{p}$, $\mathbf{y}=P\mathbf{x}$ and $x\in \ker (PA-AP)$.
\endproof

An immediate consequence is the following

\begin{corollary}\label{cor:P} Let $A$ be an $n\times n$ matrix and let  $(\mu_\mathbf{x})_{\mathbf{x}\in \mathbb{R}^n }$ be the family of measures generated by $A$. For each $P\in \mathcal{P}$ such that \begin{equation*}\label{commutatPerm}
PA=AP \quad \text{and} \quad P\mathbf{p}=\mathbf{p} ,
\end{equation*} we have that, for all $\mathbf{x}\in \R^n$,
$$\mu_\mathbf{x}=\mu_{P\mathbf{x}}.$$
\end{corollary}

We now prove a preliminary lemma that will be needed for the proof of the next theorem. \begin{lemma}\label{LemmaKern}
Fix $N\in \N \setminus \{0\}$. Given $N$ $n\times n$  non-zero matrices $K_i$, for $i\in [N]$, there exists a vector $\mathbf{v}\in \R^n$ such that $\mathbf{v}\notin \ker(K_i)$ for all $i\in [N]$.
\end{lemma}

\proof
We prove the claim by induction over $N$. For $N=1$, the claim is trivially true as the the matrix is non-zero.

We now assume the statement to be true for $N-1$, and we prove it for $N$. From the inductive hypothesis, there exists $\mathbf{v}\in \R^n$ such that $\mathbf{v}\notin \ker(K_i)$ for every $i\in [N-1]$. We also know that there exists $\mathbf{w}\in \R^n$ such that $\mathbf{w}\notin \ker(K_N)$ from the base case. Therefore, we can observe that we can choose an $\alpha >0 $ such that the vector $\mathbf{v}+\alpha \mathbf{w}\notin \ker(K_i)$ for every $i\in [N]$. In fact, for every $i\in [N]$, using linearity and the reverse triangular inequality we have
\begin{equation*}
\|K_i(\mathbf{v}+\alpha \mathbf{w})\|\geq |\|K_i\mathbf{v}\|-\alpha \|K_i\mathbf{w}\||.
\end{equation*}We now notice that $K_i\mathbf{v}$ and $K_i\mathbf{w}$ cannot both be zero for every $i\in [N]$, as a consequence of the above discussion. We can therefore choose an $\alpha$ such that the line passing trough the origin in the plane with slope $\alpha$ does not intersect any of the points with coordinates $(\|K_i\mathbf{v}\|,\|K_i\mathbf{w}\|)$.
\endproof

We also notice that in the previous lemma, we could have considered, more generally, bounded operators on a normed vector space instead of only square matrices. More generally, we could also have considered countably many bounded operators, instead of finitely many, on a complete normed vector space. In fact, the proof idea of Lemma \ref{LemmaKern} rewritten in set theoretic language reads: the kernel of a non-zero operator is nowhere dense and, therefore, the countable union of the kernels is a set of first category, therefore its complement is dense in the complete normed vector space (i.e.\ non-empty) as a consequence of Baire category theorem. \newline

The next theorem ensures the existence of irreducible vectors with pairwise distinct entries, for any given matrix.

\begin{theorem}\label{thm:Irred} Let $A$ be an $n\times n$ matrix and let  $(\mu_\mathbf{x})_{\mathbf{x}\in \mathbb{R}^n }$ be the family of measures generated by $A$. Then, 
\begin{enumerate}
    \item There exists an irreducible vector $\mathbf{v}=(v_i)_i\in \R^n$ for $A$.
\item We can always assume that the entries of $\mathbf{v}$ are pairwise distinct, that is, $v_i\neq v_j$ for $i\neq j$.
\end{enumerate}

\end{theorem}
\proof { If $PA=AP$ for every $P\in \mathcal{P}$, then every vector is an irreducible vector, and there is nothing to show. Therefore, we consider the case in which there exists at least one $P_2\in \mathcal{P}$ such that $P_2A\neq AP_2$. For every $P\in \mathcal{P}$ such that $PA\neq AP$, there exists a non-zero vector $\mathbf{v}_{P}\notin \ker (PA-AP)$. Therefore, for every $P_1\in \mathcal{P}$, $$(P_1)^{\top}\mathbf{v}_{P}\notin (P_1)^{\top}\ker (PA-AP)=\ker \bigl((PA-AP)P_1\bigr),$$where we are using the fact that $P_1^{\top}=P_1^{-1}$, since this is a permutation matrix. Now, by  Lemma \ref{LemmaKern} we can choose $\mathbf{v}\in \R^n$ such that $$\mathbf{v}\notin \ker \bigl((PA-AP)P_1\bigr)$$ for all $P\in \mathcal{P}$ such that $PA\neq AP$ and $P_1\in \mathcal{P}$, where we consider all matrices $(PA-AP)P_1$.
Therefore, given $P_1,P_2\in\mathcal{P}$, we have that $$P_1\mathbf{v}\in\ker(P_2A-AP_2) \iff P_2A=AP_2.$$
Now, by Lemma \ref{lemma:P}, we have that $\mu_{P_1\mathbf{v}}=\mu_{\mathbf{y}}$ if and only if  there exists $P_2\in \mathcal{P}$ such that $\mathbf{y}=P_2P_1\mathbf{v}$ and $P_1\mathbf{v}\in\ker(P_2A-AP_2)$. Hence, by the above observation, $\mu_{P_1\mathbf{v}}=\mu_{\mathbf{y}}$ if and only if there exists $P_2\in \mathcal{P}$ such that $\mathbf{y}=P_2P_1\mathbf{v}$ and $P_2A=AP_2$. This proves the first claim.\newline 
To prove the second claim, assume that $v_i=v_j$ for some $i,j$ such that $i\neq j$. Consider the vector
$$\mathbf{v'}\coloneqq \mathbf{v}+\xi \mathbf{e_i},$$
where the absolute value of $\xi\in\mathbf{R}$ is small enough that
$$\mathbf{v'}\notin \ker((PA-AP)P_1),$$
for $P,P_2\in \mathcal{P}$ such that $PA\neq AP$. This is possible since the kernel is closed, therefore its complement is open.
\endproof

From here on in this section, we fix a constant $K\in\mathbb{R}_{\geq 1}$ and we let $\mathcal{S}$ be the set of $n\times n$ matrices $A$ such that $||A||_{\infty \rightarrow 1}\leq K$. Moreover, given a vector $\mathbf{x}\in\mathbb{R}^n$ and $d>0$, we let 
$$\mathbf{x}^i_{d}\coloneqq \mathbf{x}+\frac{d^2}{64K}\mathbf{e_i}.$$ We observe that, if $\mathbf{x}$ satisfies the conditions of Theorem \ref{thm:Irred}, then also $\mathbf{x}^i_{d}$ satisfies them, for every $d>0$ small enough.\newline

\begin{example}
If $K\geq (n-1)$, then $\mathcal{S}$ contains all adjacency matrices associated to graphs on $n$ nodes.
\end{example}

We now recall the following result from \cite{backhausz2018action}, since it will be needed in the proof of Lemma \ref{lemmaDist1} below.

\begin{lemma}[Lemma 13.1 in \cite{backhausz2018action}]\label{lemma:backhausz}
Let $\tau(X-Y)$ denote the maximum of $\mathbb{E}\left(\left|\pi_{i}(X-Y)\right|\right)$ over $1 \leq i \leq k$, where  $\pi_{i}: \mathbb{R}^{k} \rightarrow \mathbb{R}$ is the $i$-th coordinate function for $1 \leq i \leq k$.
Let $X, Y$ be two jointly distributed $\mathbb{R}^{k}$-valued random variables. Then, \begin{equation*}
d_{\mathrm{LP}}(X_{\#}\mathbb{P},Y_{\#}\mathbb{P} ) \leq \tau(X-Y)^{1 / 2} k^{3 / 4}.
\end{equation*}\end{lemma}
We apply the above lemma to prove the following
\begin{lemma}\label{lemmaDist1}Fix $A\in\mathcal{S}$ and let $\mathbf{x}\in\mathbb{R}^n$.
Let $\nu_1\in \mathcal{Z}_{\mathbf{x}}$ and write $\nu_1=\mu_{P\mathbf{x}}$ for some $P\in\mathcal{P}$. Then, the measure $$\nu_2\coloneqq \mu_{P(\mathbf{x}_d^i)}\in \mathcal{Z}_{\mathbf{x}^i_d}$$
is such that $$\mathrm{d_{LP}(\nu_1 , \nu_2)}<\tfrac{d}{4}.$$
\end{lemma}
\proof
 We have that $$\left\|P\mathbf{x}-P(\mathbf{x}_d^i)\right\|_1\leq \frac{d^2}{64K}\left\|\mathbf{e_i}\right\|_1\leq\frac{d^2}{64K}\leq \frac{d^2}{64}$$

and $$
\left\| AP\mathbf{x}-AP(\mathbf{x}_d^i)\right\|_1\leq\left\|A\right\|_{\infty\rightarrow 1}\cdot \left\|P\mathbf{x}-P(\mathbf{x}_d^i)\right\|_{\infty}\leq\left\|A\right\|_{\infty\rightarrow 1}\cdot \frac{d^2}{64K}\leq \frac{d^2}{64}.
$$The claim then follows by letting $\mu_{P\mathbf{x}}$ and $\mu_{P(\mathbf{x}_d^i)}$ be the distributions $X_{\#}\mathbb{P}$ and $Y_{\#}\mathbb{P}$ of the $\R^2$--valued random variables $X(\omega)=(P\mathbf{x}_{\omega},AP\mathbf{x}_{\omega})$ and $Y(\omega)=(P(\mathbf{x}_d^i)_{\omega},A{P(\mathbf{x}_d^i)}_{\omega})$ from Lemma \ref{lemma:backhausz}.
\endproof

Next, we prove a theorem that allows us to associate, to any measure in $\mathcal{Z}_\mathbf{v}$, a unique measure in $\mathcal{Z}_{\mathbf{v}^i_{\varepsilon}}$, whenever $\mathbf{v}$ is irreducible and $\varepsilon$ is small enough.

\begin{theorem}\label{thm:isolMeas}
Fix $A\in\mathcal{S}$ and let $\mathbf{v}\in \R^n$ be an irreducible vector for $A$ whose entries are pairwise distinct. Let $\varepsilon>0$ be such that:
\begin{enumerate}
    \item 
$\varepsilon<\min\bigl\{\mathrm{d_{LP}}(\nu_1,\nu_2): \nu_1,\nu_2\in \mathcal{Z}_\mathbf{v} \text{ and }\nu_1\neq \nu_2\bigr\}$, where the minimum exists since $\mathcal{Z}_\mathbf{v}$ is finite;
\item $\varepsilon$ is small enough that the vectors $\mathbf{v}^i_{\varepsilon}$, for $i=1,\ldots,n$, are such that, if $P\in \mathcal{P}$ and
$$
(P\mathbf{v})_1<(P\mathbf{v})_2<\ldots<(P\mathbf{v})_n,
$$
then
$$
(P\mathbf{v}^i_{\varepsilon})_1<(P\mathbf{v}^i_{\varepsilon})_2<\ldots<(P\mathbf{v}^i_{\varepsilon})_n.
$$
\end{enumerate}
Let also $\nu_1\in \mathcal{Z}_{\mathbf{v}}$ and write $\nu_1=\mu_{P\mathbf{v}}$ for some $P\in\mathcal{P}$. Then, the measure $\nu_2\coloneqq \mu_{P(\mathbf{v}_\varepsilon^i)}$ is the unique measure $\nu\in\mathcal{Z}_{\mathbf{v}_{\varepsilon}^i}$ such that 
$$\mathrm{\mathrm{d_{LP}}(\nu_1 , \nu)}<\tfrac{\varepsilon}{4}.$$
\end{theorem}
\proof

First, we prove that $\mathcal{Z}_{\mathbf{v}}$ has cardinality larger or equal than $\mathcal{Z}_{\mathbf{v}^i_{\varepsilon}}$. Fix two distinct elements in $\mathcal{Z}_{\mathbf{v}^i_{\varepsilon}}$ and write them as $\mu_{P_1(\mathbf{v}^i_{\varepsilon})}$ and $\mu_{P_2(\mathbf{v}^i_{\varepsilon})}$, for some $P_1,P_2\in\mathcal{P}$. If we show that $\mu_{P_1\mathbf{v}}\neq \mu_{P_2\mathbf{v}}$, we are done. Assume that $\mu_{P_1\mathbf{v}}= \mu_{P_2\mathbf{v}}$. Then, since $\mathbf{v}$ is irreducible, there exists $P_3\in\mathcal{P}$ such that $P_2\mathbf{v}=P_3P_1\mathbf{v}$ and $P_3A=AP_3$. By Corollary \ref{cor:P} and by the choice of $\varepsilon$, this implies that
\begin{equation*}
  \mu_{P_1(\mathbf{v}^i_{\varepsilon})}= \mu_{P_3P_1(\mathbf{v}^i_{\varepsilon})}= \mu_{P_2(\mathbf{v}^i_{\varepsilon})},
\end{equation*}which is a contradiction since we are assuming that $\mu_{P_1(\mathbf{v}^i_{\varepsilon})}$ and $\mu_{P_2(\mathbf{v}^i_{\varepsilon})}$ are distinct. Therefore, $\mathcal{Z}_{\mathbf{v}}$ has cardinality larger or equal than $\mathcal{Z}_{\mathbf{v}^i_{\varepsilon}}$.\newline

Now, we illustrate the rest of the proof in Figure \ref{fig:distance}. By Lemma \ref{lemmaDist1} we know that $\nu_1$ and $\nu_2$ are such that $$\mathrm{d_{LP}}(\nu_1,\nu_2)<\frac{\varepsilon}{4}.$$ Fix any $\nu'_1\in \mathcal{Z}_{\mathbf{v}}$ such that $\nu'_1\neq \nu_1$, and consider the corresponding measure $\nu'_2\in\mathcal{Z}_{\mathbf{v}^i_{\varepsilon}}$ as in Lemma \ref{lemmaDist1}, so that $$\mathrm{d_{LP}}(\nu'_1,\nu'_2)<\frac{\varepsilon}{4}.$$ 
By the choice of $\varepsilon$, we have that $\mathrm{d_{LP}}(\nu_1,\nu'_1)>\varepsilon$ and now, applying the reverse triangular inequality, it is easy to see that $\mathrm{d_{LP}}(\nu_1,\nu'_2)>\frac{\varepsilon}{4}$ and $\mathrm{d_{LP}}(\nu'_1,\nu_2)>\frac{\varepsilon}{4}$.
\endproof

\begin{figure}[h]
    \centering
    \includegraphics[width=11cm]{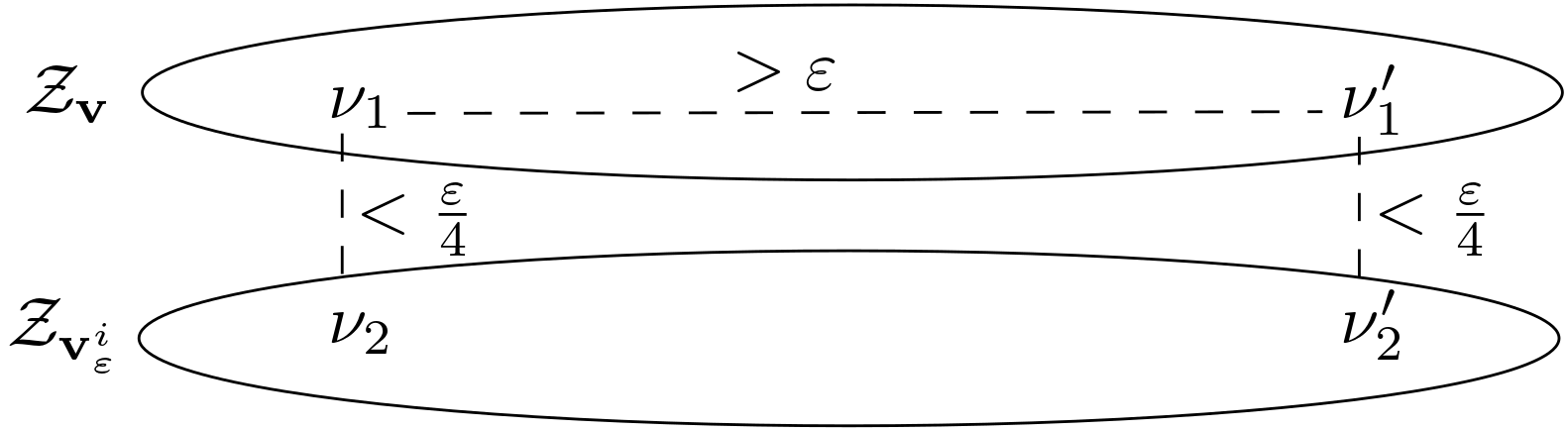}
    \caption{An illustration of the proof of Theorem \ref{thm:isolMeas}.}
    \label{fig:distance}
\end{figure}

\section{Main result }
\label{MaResults} We are now ready to prove Theorem \ref{thm:main}.

\proof[Proof of Theorem \ref{thm:main}]
If we know the set $\mathcal{Z}_A$ of measures generated by $A$, then for any $\mathbf{x}\in\mathbb{R}^n$ we also know the set $\mathcal{Z}_{(A,\,\mathbf{x})}$. Hence, we can choose a vector
 $$\mathbf{v}\in\mathrm{argmax}_{\mathbf{x}\in \R^n}\#\{\mathcal{Z}_{(A,\mathbf{x})}\}=\mathrm{argmax}_{\mathbf{x}\in \R^n}\#\{\mu_{(A,\,P\mathbf{x})}: \, P\in\mathcal{P}\}.$$
It is easy to see that $\mathbf{v}$ must be an irreducible vector for $A$ and, up to a small perturbation, we can choose $\mathbf{v}$ such that all its entries are pairwise distinct. Now, fix $\varepsilon>0$ as in Theorem \ref{thm:isolMeas}, and choose it also small enough that:
\begin{enumerate}
    \item For $i\in \{1,\ldots,n\}$, $$\mathbf{v}^i_{\varepsilon}\in\mathrm{argmax}_{\mathbf{x}\in \R^n}\#\{\mathcal{Z}_{(A,\mathbf{x})}\},$$ so that also these vectors are irreducible for $A$;
\item The entries of each $\mathbf{v}^i_{\varepsilon}$ are pairwise distinct. 
\end{enumerate}

Fix now $\nu_1\in \mathcal{Z}_{(A,\,\mathbf{v})}$. Then, we know that $\nu_1=\mu_{(A,\,P_1\mathbf{v})}$, for some $P_1\in\mathcal{P}$, but we don't know $A$ and $P_1$. However, we know that there exists $P_2\in\mathcal{P}$ such that the unique ordered support of $\nu_1$ is (cf.\ Example \ref{ex:proof} below)
\begin{equation}\label{represMeas1}(P_2P_1\mathbf{v},\,P_2AP_1\mathbf{v}),\end{equation}and we do know the resulting pair in \eqref{represMeas1}. In particular, since we know both $\mathbf{v}$ and $P_2P_1\mathbf{v}$, we can also reconstruct $P_2P_1$.\newline

Now, fix $i\in \{1,\ldots,n\}$ and observe that, since $P_1,P_2\in\mathcal{P}$, there exists $j\in \{1,\ldots,n\}$ such that $\mathbf{e_j}=P^\top_1P_2^\top\mathbf{e_i}$. By applying Theorem \ref{thm:isolMeas} to $\nu_1$, we know that the measure  $$\nu_2\coloneqq \mu_{\left(A,\,P_1(\mathbf{v}^j_{\varepsilon})\right)}\in \mathcal{Z}_{\left(A,\,\mathbf{v}^j_{\varepsilon}\right)}$$
is the unique measure $\nu\in \mathcal{Z}_{(A\,,\mathbf{v}^j_{\varepsilon})}$ such that 
$$\mathrm{\mathrm{d_{LP}}(\nu_1 , \nu)}<\tfrac{\varepsilon}{4}.$$

Hence, we can identify $\nu_2$ and therefore also its unique ordered support, that we can write as \begin{align}\label{repMeas2}
&\left(P_2P_1(\mathbf{v}^j_{\varepsilon}),\,P_2 AP_1(\mathbf{v}^j_{\varepsilon}))\right) \nonumber \\
&=\left(P_2P_1(\mathbf{v}+\frac{\varepsilon^2}{64K} \mathbf{e_j}),\,P_2 AP_1(\mathbf{v}+\frac{\varepsilon^2}{64K} \mathbf{e_j})\right) \nonumber \\
&=\left(P_2P_1\mathbf{v}+\frac{\varepsilon^2}{64K} \mathbf{e_i},\,P_2AP_1\mathbf{v}+\frac{\varepsilon^2}{64K} P_2AP^\top_2\mathbf{e_i}\right).
\end{align}Taking the difference between \eqref{represMeas1} and \eqref{repMeas2} leads to
\begin{equation}
\frac{\varepsilon^2}{64K}(\mathbf{e_i},P_2AP^\top_2\mathbf{e_i}),
\end{equation}
 therefore we are able to reconstruct the $i-$th column of $P_2 AP_2^\top$. Since we can do this for every $i\in\{1,\ldots,n\}$, we can reconstruct the entire matrix $P_2AP^\top_2$. This proves the claim.
\endproof

We illustrate the first part of the proof of Theorem \ref{thm:main} with an example.
\begin{example}\label{ex:proof}
Let $n=3$, so that $|\mathcal{P}|=6$, and let
\begin{equation*}
    A\coloneqq\begin{pmatrix}
    0 & 1 & 1\\
    0 & 1 & 0 \\
    1 & 1 & 0
    \end{pmatrix}, \quad P'\coloneqq\begin{pmatrix}
    0 & 0 & 1\\
    0 & 1 & 0 \\
    1 & 0 & 0
    \end{pmatrix}\in \mathcal{P}.
\end{equation*}Then, $P'A=AP'$, implying that
\begin{equation*}
    \max_{\mathbf{x}\in \R^n}\#\{\mathcal{Z}_{(A,\mathbf{x})}\}=\max_{\mathbf{x}\in \R^n}\#\{\mu_{(A,\,P\mathbf{x})}: \, P\in\mathcal{P}\}\leq 3.
\end{equation*}Now, by letting
\begin{equation*}
\mathbf{v}\coloneqq \begin{pmatrix}v_1 \\ v_2 \\ v_3\end{pmatrix}\coloneqq \begin{pmatrix}3 \\ 1 \\ 2\end{pmatrix},  
\end{equation*}it is easy to see that
\begin{equation*}
    \#\{\mathcal{Z}_{(A,\mathbf{v})}\}=\#\{\mu_{(A,\,P\mathbf{v})}: \, P\in\mathcal{P}\}= 3,
\end{equation*}therefore
$$\mathbf{v}\in\mathrm{argmax}_{\mathbf{x}\in \R^n}\#\{\mathcal{Z}_{(A,\mathbf{x})}\}=\mathrm{argmax}_{\mathbf{x}\in \R^n}\#\{\mu_{(A,\,P\mathbf{x})}: \, P\in\mathcal{P}\},$$
hence, in particular, $\mathbf{v}$ is irreducible for $A$. The support of $\mu_{(A,\,\mathbf{v})}$ is illustrated in Figure \ref{fig:example}(a). Now, let
\begin{equation*}
    P_1\coloneqq\begin{pmatrix}
    0 & 1 & 0\\
    1 & 0 & 0 \\
    0 & 0 & 1
    \end{pmatrix}\in\mathcal{P}.
\end{equation*}

\begin{figure}[h]
    \centering
    \subfloat[\centering ]{{\includegraphics[width=5cm]{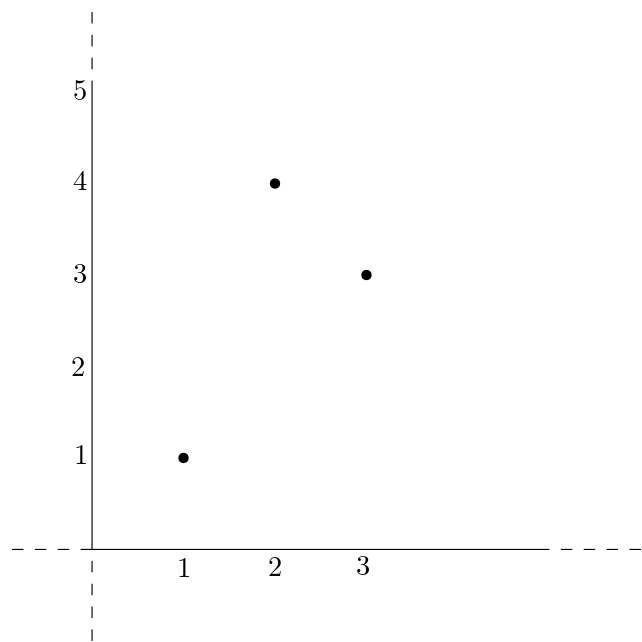} }}
    \qquad
    \subfloat[\centering ]{{\includegraphics[width=5cm]{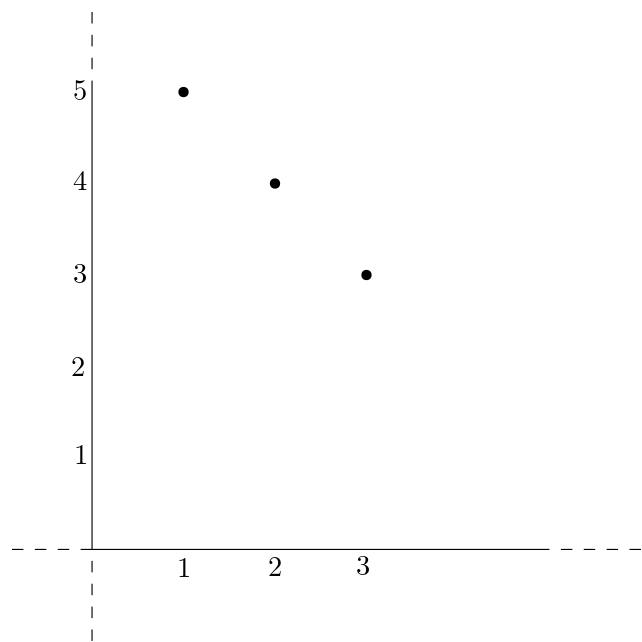} }}
    \caption{The supports of $\mu_{(A,\,\mathbf{v})}$ and $\mu_{(A,\,P_1\mathbf{v})}$ in Example \ref{ex:proof}.}
    \label{fig:example}
\end{figure}

Figure \ref{fig:example}(b) shows the support of $\nu_1:=\mu_{(A,\,P_1\mathbf{v})}$, and it is clear from the picture that the unique ordered support of $\nu_1$ is
\begin{equation*}
\left(\begin{pmatrix}1 \\ 2 \\ 3\end{pmatrix},\,\begin{pmatrix}5 \\ 4 \\ 3\end{pmatrix} \right)=(P_2P_1\mathbf{v},\,P_2AP_1\mathbf{v}),
\end{equation*}where
\begin{equation*}
    P_2\coloneqq\begin{pmatrix}
    1 & 0 & 0\\
    0 & 0 & 1 \\
    0 & 1 & 0
    \end{pmatrix}\in\mathcal{P}.
\end{equation*}
\end{example}
We note that, if we consider measures $\mu_x$ which are relative to a probability measure $\mathbf{p}$ that is not the uniform probability measure, then a matrix $A$ as in Theorem \ref{thm:main} will be characterized by a stronger notion than switching equivalence. In particular, in this case $\mathcal{Z}_A=\mathcal{Z}_B$ if and only if $A=PBP^{\top}$, where $P$ is such that $P\mathbf{p}=\mathbf{p}$. An interesting particular case is when the vector $\mathbf{p}$ has pairwise distinct entries $p_i\neq p_j$ for $i\neq j$. In this case, $\mathcal{Z}_A=\mathcal{Z}_B$ if and only if $A=B$.
\section{Properties of matrices and graphs from the generated measure }
\label{PropertiesGraphs}
We now discuss how some properties of matrices and graphs directly translate in terms of the measure generated by the matrices. \newline

We first notice that, for an $n\times n$ matrix $A$ and a scalar $\lambda\in \mathbb{R}$, a measure in $\mathcal{Z}_A$ is supported on the graph of the linear function $\lambda x$, $$\mathrm{graph}(\lambda x)\coloneqq\{(x,\lambda x)\in \mathbb{R}^2: \ \lambda\in \mathbb{R}\},$$if and only if $\lambda$ is an eigenvalue of the matrix $A$.
Moreover, we notice that the measure $\mu$ which is completely supported on the vertical line $$\{(1, x)\in \R^2: \ x\in \R \}$$ corresponds to the constant $1$ vector. For this reason, the set $$\{x\in \R : \  (1, x) \in \mathrm{supp}(\mu)\}$$ corresponds to the row sums of the row of the matrix $A$.\newline

We now want to relate our measure representation of matrices to graphs. In order to do this, we introduce several matrix representations of graphs. We start with the simplest possible choice.
\begin{definition}
Let $G=(V,E)$ a graph and on vertices $v_{1}, \ldots, v_{N}$. The \emph{adjacency matrix} of $G$ is the $N\times N$ matrix $A$, whose entries are
$$
A_{i j}:= \begin{cases}1 & \text { if }\left(v_{i}, v_{j}\right) \in E \\ 0 & \text { otherwise }.\end{cases}
$$
\end{definition} 
Another possible matrix representation that can be advantageous in some cases is the following.
\begin{definition}
The \emph{Kirchhoff Laplacian} of $G$ is the $N\times N$ matrix
$$
K:=D-A,
$$
where $D$ is the diagonal matrix of the degrees.
\end{definition}

We remark that both the adjacency matrix and the Kirchhoff Laplacian matrix are such that two graphs are isomorphic if and only if the respective matrices $M_1$ and $M_2$ are related by the the relationship $M_1=PM_2P$.\newline

Additionally, we recall the following

\begin{definition}
The \emph{normalized Laplacian} is the $N\times N$ matrix 
$$
L:=\mathrm{Id}-D^{-1 } A, 
$$ 
where $\mathrm{Id}$ is the $N \times N$ identity matrix. For $i\neq j$, we have
$$
L_{i j}=-\frac{A_{i j}}{\operatorname{deg} v_{i}},
$$
which is minus the probability that a random walker goes from $v_{i}$ to $v_{j}$.
\end{definition}

We refer to \cite{chung1997spectral,MHJ} for more details on the normalized Laplacian and its spectral theory.\newline

We now consider the above matrices in relation to graphs, and we observe how properties of graphs translate in the family of measure representations. \newline

We consider a graph $G=(V,E)$ and the related adjacency matrix $A$. We already discussed how to directly extract the spectral information of a matrix. Additionally, as the degrees of the graph $G$ correspond to the sums of the rows of the adjacency matrix, we can directly determine from $\mathcal{Z}_A$ the degrees as presented above.\newline

The spectrum and the degree distribution determine many graph properties. For example, the spectrum of the normalized Laplacian determines the number of connected components, and whether the graph is bipartite. Additional properties are characterized by this information. We now recall the following definition, and we refer to \cite{HomomChayes} for more details on this topic. 

\begin{definition}
For two graphs $F=(V(F),E(F))$ and $G=(V(G),E(G))$, a \emph{graph homomorphism} from $F$ to $G$ is a function $\phi: V(F) \rightarrow V(G)$ such that for $v,w\in V(F)$ and $\{v,w\}\in E(F)$ implies $\{\phi(u),\phi(v)\}\in E(G)$. We denote by $\operatorname{hom}(F,G)$ the number of homomorphisms from $F$ to $G$.
\end{definition}

We now give some examples of homomorphism numbers that can be extracted directly from the measure--theoretic representation. 

\begin{example}Let $S_{k}$ denote the star graph with $k$ nodes. Then, for any graph $G$ on $n$ nodes,
$$
\operatorname{hom}\left(S_{k}, G\right)=\sum_{i=1}^{n} d_{i}^{k-1},
$$
where $d_{1}, \ldots, d_{n}$ are the degrees of $G$. 
\end{example}

\begin{example}
Let $C_{k}$ denote the cycle graph on $k$ nodes, and again let $G$ be any graph on $n$ nodes. Then 
$$
\operatorname{hom}\left(C_{k}, G\right)=\sum_{i=1}^{n} \lambda_{i}^{k}
$$
where $\lambda_{1}, \ldots, \lambda_{n}$ are the eigenvalues of the adjacency matrix of $G$. 
\end{example}

Notice that, in general, it is not easy to reconstruct directly general information about the adjacency matrix $A$ and its powers, because the measure representation does not keep track of permutations of vectors.

\section{Relationship with action convergence metric }
\label{RelatActionConv}

In this section we briefly present the notion of action convergence from \cite{backhausz2018action}, and we underline the relationship with the metric defined in this work. We start by recalling the following
\begin{definition}
 A \emph{P-operator} is a linear operator of the form $A: L^{\infty}(\Omega) \rightarrow L^{1}(\Omega)$ such that
$$
\|A\|_{\infty \rightarrow 1}:=\sup _{X \in L^{\infty}(\Omega), \substack{ X\neq 0}}\frac{\|AX\|_{1}}{\|X\|_{\infty}} 
$$
is finite, where $(\Omega,\mathcal{F}, \mathbb{P})$ is a generic probability space. We denote by $\mathcal{B}(\Omega)$ the set of all $P$-operators on $\Omega$.
\end{definition}

We show, with the following example, that a matrix is a $P-$operator.
\begin{example}
 If $\Omega=[n]:=\{1,\ldots,n\}$ and $\mathbb{P}=(p_{\omega})_{\omega\in[n]}$ on $\Omega$ is the probability measure relative to the probability vector $\mathbf{p}=(p_{\omega})_{\omega\in[n]}$, then $L^{1}(\Omega)=L^{\infty}(\Omega)=\mathbb{R}^{n}$.
In this case, $\mathcal{B}(\Omega)$ is the set of all $n \times n$ matrices. Thus, every matrix $A \in \mathbb{R}^{n \times n}$ is a $P$-operator.
\end{example}
We consider now a $P-$operator
$$
A: L^{\infty}(\Omega)\rightarrow L^1(\Omega)
$$
and we let
$$
X\in L^{\infty}(\Omega).
$$
This is a bounded random variable, and we also have that
$$
AX\in L^1(\Omega)
$$
is a random variable with finite expectation. \newline
We can therefore define the $2-$dimensional random vector
$(X,AX).$\newline
More generally, for all random variables $Z_1,Z_2,\ldots,Z_k\in L^{\infty}(\Omega)$, we can define the $2k-$dimensional random vector
$$
(Z_1,AZ_1,Z_2,AZ_2,\ldots,Z_k,AZ_k).
$$

We show, with the following example, how the $2-$dimensional random vectors constructed above relate to the measures generated by $A$ defined in Section \ref{Definitions}.
\begin{example} For the probability space $\Omega=[n]$, with probability $\mathbb{P}=(p_{\omega})_{\omega\in[n]}$, a matrix $A\in \mathbb{R}^{n\times n}$ and a vector $X=(X(\omega))_{\omega\in [n]}\in \mathbb{R}^n$ the law of the $2-$dimensional random vector $(X,AX)$ is 
$$
\sum_{\omega\in [n]}p_{\omega} \delta_{\left(X(\omega),(AX)(\omega)\right)}.
$$These, in fact, are the measures generated by $A$ according to Definition \ref{defMeasGenerated}.

\end{example}

We now define the $k-$profile of $A$ as 
$$S_k(A)=\bigcup_{Z_1,\ldots,Z_k\in L_{[-1,1]}^{\infty}(\Omega)}\{((Z_1,AZ_1,Z_2,AZ_2,\ldots,Z_k,AZ_k))_{\#}\mathbb{P}\}.
$$ 
Considering two $P-$operators
$$
A: L^{\infty}(\Omega_1)\rightarrow L^1(\Omega_1)
$$and
$$
B: L^{\infty}(\Omega_2)\rightarrow L^1(\Omega_2),
$$
we finally define the action convergence metric.
\begin{definition}[Metrization of action convergence] For the two $P$-operators $A, B$ the \emph{action convergence distance} is
$$
d_{M}(A, B):=\sum_{k=1}^{\infty} 2^{-k} d_{H}\left(\mathcal{S}_{k}(A), \mathcal{S}_{k}(B)\right)
$$
\end{definition}

This metric has some nice compactness properties as stated in the following
\begin{theorem}[Theorem 2.9 in \cite{backhausz2018action}]
 Let $p \in[1, \infty)$ and $q \in[1, \infty]$. Let $\left\{A_{i}\right\}_{i=1}^{\infty}$ be a convergent sequence of $P$-operators with uniformly bounded $\|\cdot\|_{p \rightarrow q}$ norms. Then there is a $P$-operator A such that $\lim _{i \rightarrow \infty} d_{M}\left(A_{i}, A\right)=0$, and $\|A\|_{p \rightarrow q} \leq \sup _{i \in \mathbb{N}}\left\|A_{i}\right\|_{p \rightarrow q}$.
\end{theorem}

Moreover, action convergence unifies several approaches to graph limits theory. In particular, consider the sequence of adjacency matrices $A_n$ of graphs $G_n$, and let $v_n$ be the number of vertices of $G_n$. Then,  \begin{itemize}
    \item The action convergence of the sequence $$
   \frac{A_n}{v_n} 
    $$
    coincides with graphon convergence \cite[Theorem 8.2 and Lemma 8.3]{backhausz2018action}
    \item The action convergence of the sequence $$
    A_n
    $$
    coincides with local-global convergence \cite[Theorem 9.2]{backhausz2018action}.\end{itemize}

We have that the following special class of $P-$operators is important in graph limits theory as it naturally generalizes the notion of adjacency matrix of a graph. 
\begin{definition}
A positivity-preserving and self-adjoint $P-$operator
is called a \emph{graphop}.
\end{definition}
Consider two $P-$operators
$$
A: L^{\infty}(\Omega_1)\rightarrow L^1(\Omega_1)
$$and
$$
B: L^{\infty}(\Omega_2)\rightarrow L^1(\Omega_2).
$$We can now define the simplified metric
\begin{definition}[$1-$profile Metric] For the two $P$-operators $A, B$ the \emph{1-profile distance} is
$$
d_{S}(A, B):= d_{H}\left(\mathcal{S}_{1}(A), \mathcal{S}_{1}(B)\right)
$$
\end{definition}
We now compare the $1-$profile Metric introduced in this work with the action convergence metric. We have the following

\begin{lemma}
$$d_{S}(A, B)\leq 2d_{M}(A, B)$$
\end{lemma}
\proof
\begin{equation*}
\begin{aligned}
d_{S}(A, B)&=d_{H}\left(\mathcal{S}_{1}(A), \mathcal{S}_{1}(B)\right)\\
&\leq d_{H}\left(\mathcal{S}_{1}(A), \mathcal{S}_{1}(B)\right)+2\sum_{k=2}^{\infty} 2^{-k} d_{H}\left(\mathcal{S}_{k}(A), \mathcal{S}_{k}(B)\right)\\
& \leq 2d_{M}(A, B).
\end{aligned}
\end{equation*}
\endproof

Moreover, a direct Corollary of Theorem \ref{MaResults} is the following
\begin{corollary}
Let $K\in \mathbb{R}_{\geq 1}$. Let $\mathcal{S}$ be the set of $n\times n$ matrices $A$ such that $||A||_{\infty \rightarrow 1}\leq K$. Then, for matrices $A, B\in \mathcal{S}$ we get that 

\begin{equation*}
d_{M}(A, B)=0
\end{equation*}
if and only if 
\begin{equation*}
d_{L}(A, B)=0
\end{equation*}
\end{corollary}
 Notice that, on the space of graph isomorphisms of finite graphs, the $1-$profile metric and the action convergence metric both induce the discrete topology and, therefore, we have the following 
\begin{corollary}
The $1-$profile metric and the action convergence metric are topologically equivalent on the space of graph isomorphisms for graphs with at most $n$ vertices.
\end{corollary}
This is obviously not clear anymore for general $P-$operators as the topology induced by the two metrics is not anymore discrete.
\section{Future directions}
\label{Concl}

In future work, we aim to better understand when the action convergence metric is topologically equivalent to the simplified $1$-profile Metric we introduced. This would contribute to a better understanding of action convergence, potentially giving new insight about the difference between the convergence of dense graph sequences and sparse/bounded degree sequences. The $1$-profile metric could be related to weaker notions of local-global convergence for bounded degree graph sequences, as the notion of Benjamini-Schramm convergence. A better understanding of these types of metrics could also help to understand if limiting the number of colorings in the notion of local-global convergence to 2 or more does actually change the notion of convergence, as asked in \cite{Hatami2014LimitsOL}.

\section*{Funding}
Raffaella Mulas was supported by the Max Planck Society's Minerva Grant. 

\section*{Acknowledgments}
The authors would like to thank Tobias Böhle, Christian Kuehn, Florentin M\"unch and Jiaxi Nie for helpful discussions. They are grateful to the anonymous referee for the comments and suggestions that have greatly improved the first version of this paper.

\bibliographystyle{plain} 
\bibliography{Bibliography}

\end{document}